\documentclass[12pt,a4paper,reqno]{amsart}

\usepackage{tcolorbox}

\usepackage{amsmath,amsfonts,amssymb,latexsym,graphics,epsfig,url}
\usepackage[dvipsnames]{xcolor}
\usepackage[normalem]{ulem}
\usepackage{amsthm}
\usepackage[english]{babel}
\usepackage{graphicx,tikz}
\usetikzlibrary{decorations.markings}
\usepackage{soul}
\usepackage[utf8]{inputenc}
\usepackage[T1]{fontenc}
\usepackage{hyperref} 
\hypersetup{colorlinks,
    linkcolor=blue, 
    anchorcolor=blue,
    citecolor=blue}
\tikzstyle{directed}=[postaction={decorate,
	decoration={markings,mark=at position 0.6 with {\arrow{latex}}}
}]
\tikzstyle{vertex}=[circle, draw, inner sep=0pt, minimum size=3pt]

\usepackage{algpseudocode}
\usepackage{amsthm}
\usepackage{fullpage}
\usepackage{float}
\usepackage{multirow}
\usepackage{amscd}
\usepackage{array}
\usepackage{hyperref}
\usepackage{mleftright} 
\usepackage{float}
\usepackage{ulem}
\usepackage{soul}
\PassOptionsToPackage{lowtilde}{url}
\usepackage{xcolor}
\usepackage{cancel}
\usepackage{tikz}

\newtheorem{theorem}{Theorem}[section]
\newtheorem{lemma}[theorem]{Lemma}
\newtheorem{corollary}[theorem]{Corollary}
\newtheorem{proposition}[theorem]{Proposition}

\newtheorem{remark}[theorem]{Remark}

\def\1{{\bf 1}}

\def\vu{\mbox{\boldmath $u$}}
\def\vv{\mbox{\boldmath $v$}}

\def\Cay{\mbox{\rm Cay}}

\DeclareMathOperator{\diam}{diam}
\DeclareMathOperator{\spec}{sp}

\def\BB{\mbox{\boldmath $B$}}
\def\vv{{\mbox {\boldmath $v$}}}

\def\Z{\mathbb Z}
\newcommand{\cal}{\mathcal}

\def\blue{\textcolor{blue}}

\def\purple{\textcolor{purple}}

\def\greeen{\textcolor[RGB]{50, 220, 26 }}

\date{}
\begin{document}
\title{On generalized token graphs with tokens of different colors}

\author{Cristina Dalf\'o}
\address{Departament de Matem\`atica, Universitat de Lleida, Igualada (Barcelona), Catalonia}
\email{cristina.dalfo@udl.cat}

\author{Blas Fern\'andez}
\address{University of Primorska, Koper, Slovenia}
\email{blas.fernandez@famnit.upr.si}

\author{Miquel \`Angel Fiol}
\address{Departament de Matem\`atiques, Universitat Politècnica de Catalunya, Barcelona, Catalonia; Barcelona Graduate School of Mathematics, Barcelona, Catalonia; Institut de Matem\`atiques de la UPC-BarcelonaTech (IMTech), Barcelona, Catalonia}
\email{miguel.angel.fiol@upc.edu}

\author{M\'onica A. Reyes}
\address{Departament de Matem\`atica, Universitat de Lleida, Igualada (Barcelona), Catalonia}
\email{monicaandrea.reyes@udl.cat}

\begin{abstract}
We generalize the concept of token graphs, and we call them supertoken graphs. In this work, we consider the variant in which tokens are distinguishable (for instance, by color) and at most one token is allowed on each vertex. We formally define this kind of supertoken graphs and present their basic properties. We study their diameter, in which cases they are Cayley graphs, and their spectra. Moreover, we 
provide bounds or exact values for the independence, the vertex cover, the clique, and the chromatic numbers of these graphs. Finally, we characterize the cases in which the supertoken graphs of complete graphs are planar.
\end{abstract}

\maketitle



\section{Introduction}
\label{sec:supertoken}

Given a graph $G$ on $n=|V(G)|$ vertices and an integer $k\ge 1$, the \textit{$k$-supertoken graph} $F^{1}_{k\times 1}(G)$ is the graph whose vertices are the $k$-tuples $\vu=(u_1,u_2,\ldots,u_k)$ (also denoted $u_1u_2\ldots u_k$) such that $u_i\in V(G)$ for all $i$ and $u_i\neq u_j$ whenever $i\neq j$. Thus, each vertex of $F^{1}_{k\times 1}(G)$ represents a placement of $k$ distinguishable tokens on pairwise distinct vertices of $G$.
Two vertices $\vu=(u_1,\ldots,u_k)$ and $\vv=(v_1,\ldots,v_k)$ are adjacent if there exists a unique index $i$ such that $u_j=v_j$ for every $j\neq i$ and $u_iv_i\in E(G)$. In other words, $\vv$ is obtained from $\vu$ by moving exactly one token along an edge of $G$ to a vertex not occupied by any other token. Hence, $\vu$ and $\vv$ have exactly $k-1$ common entries.



Supertoken graphs are a particular case of a broader family that we call 
\textit{generalized token graphs}.
In such graphs, each vertex represents a configuration of $k$ (undistinguished or distinguished) tokens placed on some (different or equal) vertices of $G$. Moreover, two vertices (or configurations) are adjacent when one can be obtained from the other by moving one token in $G$ along an edge. Depending on the nature and positions of the tokens and the adjacency rules, we obtain the following families of generalized token graphs.

\begin{itemize}
\item
$F_k(G)$: The \textit{$k$-token graph} of $G$ (Fabila-Monroy, Flores-Pe\~naloza, Huemer, Hurtado, Urrutia, and Wood \cite{ffhhuw12}) or \textit{$k$-th symmetric power} of $G$ (Audenaert, Godsil, Royle, and Rudolph \cite{agrr07}). It has at most one token per vertex, and the tokens are indistinguishable.
\item
${\cal F}_k(G)$: The \textit{$k$-supertoken} graph of $G$ (Baskoro, Dalf\'o, Fiol, and Simanjuntak \cite{bdfs24}) or \textit{$k$-th reduced power graph} of $G$ (Hammack and Smith \cite{HaSm17}). It has at most $k$ tokens per vertex, and the tokens are indistinguishable.
\item 
$F^{1}_{k\times 1}(G)$: The \textit{$k$-supertoken graph} of $G$ with tokens of different colors (Song, Dalfó, Fiol, Mora, and Zhang in \cite{sdfmz26}), with at most one token per vertex.
\item 
$F^{s}_{k\times 1}(G)$: The \textit{$k$-supertoken graph} of $G$ with tokens of different colors (Song, Dalfó, Fiol, Mora, and Zhang in \cite{sdfmz26}). It has at most $s(<k)$ tokens per vertex.
\item 
The \textit{$k$-supertoken graph} $F^{k}_{k\times 1}(G)$ (Song, Dalfó, Fiol, Mora, and Zhang in \cite{sdfmz26}) has at most $k$ tokens per vertex, and the tokens are distinguishable. $F^{k}_{k\times 1}(G)$ is isomorphic to the Cartesian product $F^{k}_{k\times 1}(G)\cong G\Box\stackrel{(k)}{\cdots}\Box G$.
\end{itemize}

This paper is structured as follows. In the following section, we present the known results for this type of supertoken graphs. 
In Section \ref{sec:prop}, we present some fundamental properties of supertoken graphs. 
We study their diameter in Section \ref{sec:diam}, in which cases they are Cayley graphs in Section \ref{sec:Cayley}, and their spectra in Section \ref{sec:spectra}. Moreover, we provide bounds or exact values for the independence and the vertex cover numbers, both in Section \ref{sec:ind-number}. In Sections \ref{sec:clique} and \ref{sec:chrom}, we give the clique and the chromatic numbers of these graphs, respectively. Finally, in the last section, we characterize the cases in which the supertoken graphs of complete graphs are planar.


\section{Known results of supertoken graphs}

In this type of supertoken graph, the $k$ tokens are distinct, and at most one token is allowed on each vertex. 
Thus, the vertices of $F^{1}_{k\times 1}(G)$ correspond to ordered $k$-tuples of distinct 
vertices of $G$. Consequently, if $G$ has order $n$ and size $m$, then the supertoken graph 
$F^{1}_{k\times 1}(G)$ has order $\frac{n!}{(n-k)!}$ and size $mk\frac{(n-2)!}{(n-k-1)!}$. 

Some examples of this construction are shown in Figures~\ref{fig:F_2x1^1(K5)}, \ref{fig:F_2x1^1(S3)}, \ref{fig:F_3x1^1(K4)}, \ref{fig:F_2(S4)}(right), \ref{fig:F_(2x1)^1(C_7)}, and \ref{fig:F21K4}.
Some basic families of supertoken graphs have already been described in the literature. 
In particular, the supertoken graphs of paths and cycles were studied by Song, Dalfó, 
Fiol, Mora, and Zhang \cite{sdfmz26}. We recall the corresponding results below, which 
describe when cycles appear in these graphs.



\begin{proposition}[\cite{sdfmz26}]
\label{prop:path-cycle}
    Let $G=P_n$ be a path with order $n$. Then, $F^{1}_{k\times 1}(G)$ contains at least one cycle for $k\le n-2$, and does not contain any cycle for $k=n-1,n$. Moreover, for $k=n$, $F^{1}_{k\times 1}(G)$ is the graph with $n!$ isolated vertices, and $F^{1}_{k\times 1}(G)=P_n \cup \stackrel{((n-1)!)}{\cdots}\cup P_n$ for $k=n-1$.
\end{proposition}

For example, the supertoken graph $F^{1}_{3\times 1}(P_4)$ is $P_4\cup P_4\cup P_4\cup P_4\cup P_4\cup P_4$.

A similar description holds when the original graph is a cycle. 
The next result provides the structure of $F^{1}_{k\times 1}(C_n)$ in this case.

\begin{proposition}[\cite{sdfmz26}]
\label{prop:cycle-cycle}
    Let $G=C_n$ be a cycle of order $n$. Then $F^{1}_{k\times 1}(G)$ contains at least one cycle for $k\le n-1$, and does not contain any cycle for $k=n$. Moreover, for $k=n$, the supertoken graph $F^{1}_{k\times 1}(G)$ is the graph with $n!$ isolated vertices and $F^{1}_{k\times 1}(G)=C_{n(n-1)} \cup \stackrel{((n-2)!)}{\cdots}\cup C_{n(n-1)}$ for $k=n-1$.
\end{proposition}

For instance, the supertoken graph $F^{1}_{3\times 1}(C_4)$ is $C_{12}\cup C_{12}$.

\section{Fundamental properties of supertoken graphs}
\label{sec:prop}

We begin by establishing some basic structural properties of supertoken graphs. 
The first remark shows that, in an extreme case, the construction becomes independent 
of the original graph. We then prove that bipartiteness is inherited by the 
corresponding supertoken graph, although the converse does not hold in general.


\begin{remark}
\label{prop:any-graph}
    Let $G$ be a graph with order $n$. For $k=n$, the $F^{1}_{k\times 1}(G)$ is the graph with $n!$ isolated vertices, independently of the original graph $G$.
\end{remark}


\begin{proposition}
\label{prop:bip}
Let $G$ be a bipartite graph. Then, its $k$-supertoken graph $F^{1}_{k\times 1}(G)$ is also bipartite. 
\end{proposition}
\begin{proof}
The proof follows the same lines of reasoning as in Fabila-Monroy, Flores-Peñaloza, Huemer, Hurtado, Urrutia, and Wood \cite{ffhhuw12}, but now using vectors instead of sets.
Let $A,B$ be vectors (representing vertices of $F^{1}_{k\times 1}(G)$) with elements in $C_1\cup C_2$, where $C_1$ and $C_2$ are the stable sets of $G$. Then, the stable sets of the $k$-supertoken graph $F^{1}_{k\times 1}(G)$ are 
\begin{align*}
{\cal S}_1 &=\{A : |A|=k,  \mbox{ $A$ has an even number of elements of $C_1$}\},\\
{\cal S}_2 &=\{B : |B|=k,  \mbox{ $B$ has an odd number of elements of $C_1$}\}.
\end{align*}
Indeed, since $G$ is bipartite, the vertices adjacent to a vertex $A\in {\cal S}_1$ are obtained by `moving' a token from $C_1$ to $C_2$ {or from $C_2$ to $C_1$. In either case,} this changes the parity of the number of tokens in $C_1$ and, hence, such vertices must belong to ${\cal S}_2$. Analogously, all the vertices adjacent to vertex $B\in {\cal S}_2$ must be in
${\cal S}_1$.
\end{proof}
The converse of Proposition \ref{prop:bip} does not hold. For instance, it can be shown that the 2-supertoken $F_{2\times 1}^1(C_{2n+1})$ of the (non-bipartite) odd cycle $C_{2n+1}$ is bipartite. See Figure \ref{fig:F_(2x1)^1(C_7)} for the case of $C_7$. In fact, when $G=K_n$, we have the following characterization.

\begin{proposition}
Let \(1\le k\le n\). Then the graph \(F^{1}_{k\times1}(K_n)\) is bipartite if and only if
\[
k\in\{n-1,n\}.
\]
\end{proposition}
\begin{proof}
We consider the graph \(F^{1}_{k\times1}(K_n)\) in the distinguishable-token setting, with no two tokens allowed to occupy the same vertex.
First, suppose that \(1\le k\le n-2\). We show that \(F^{1}_{k\times1}(K_n)\) is not bipartite by exhibiting a triangle.
Fix \(k-1\) distinct vertices
\[
a_1,a_2,\dots,a_{k-1}\in V(K_n).
\]
Since \(k\le n-2\), there remain at least three further distinct vertices, say
\[
b_1,b_2,b_3\in V(K_n)\setminus\{a_1,\dots,a_{k-1}\}.
\]
Now consider the three vertices of \(F^{1}_{k\times1}(K_n)\)
\[
u=(a_1,\dots,a_{k-1},b_1),\quad
v=(a_1,\dots,a_{k-1},b_2),\quad
w=(a_1,\dots,a_{k-1},b_3).
\]
Any two of these vertices differ in exactly one coordinate. Since the original graph is complete, moving the last token from \(b_i\) to \(b_j\) is always allowed. Hence, \(u\), \(v\), and \(w\) are pairwise adjacent, and therefore they induce a triangle. Thus \(F^{1}_{k\times1}(K_n)\) is not bipartite.
Next, suppose that \(k=n\). In this case, there is no empty vertex, so no token can move. Hence, \(F^{1}_{n\times1}(K_n)\) consists of \(n!\) isolated vertices, and is therefore bipartite.
Finally, the case \(k=n-1\) is considered in Section \ref{sec:Cayley} (see Proposition \ref{propo51}).
This proves that \(F^{1}_{k\times1}(K_n)\) is bipartite if and only if \(k=n-1\) or \(k=n\).
\end{proof}

\section{The diameter of supertoken graphs}
\label{sec:diam}
As it is well known, the token graph $F_k(K_n)$ is isomorphic to the Johnson graph $J(n,k)$. Thus, from Proposition \ref{prop:regularF1_1111} (see Section \ref{sec:spectra}),
$F_{k\times 1}^1(K_n)$ admits a regular partition with quotient graph $J(n,k)$. If $k\le n/2$ (we can suppose this since $J(n,k)\cong J(n,n-k)$) the diameter of $J(n,k)$ is $k$. Then, as expected, the diameter of $F_{k\times 1}^1(K_n)$ is greater than this, as shown in the following result.
\begin{proposition}
\label{propo:diamF2(Kn)}
The diameter of the $k$-supertoken graph of $G=K_n$ is
$$
\diam (F_{k\times 1}^1(K_n))=\left\lfloor \frac{3k}{2} \right\rfloor=
\left\{
\begin{array}{ll}
3r & \mbox{if $k=2r$,}\\
3r+1 & \mbox{if $k=2r+1$.}
\end{array}
\right.
$$
\end{proposition}
\begin{proof}
Since $F_{k\times 1}^1(K_n)$ is vertex-transitive, it is enough to consider the distances from a generic vertex $x_1x_2\ldots x_k$, with $x_i\in [n]:=\{1,2,\ldots,n\}$,
to vertex $12\ldots k$. If, for some $i$, with $x_j\neq i$ and   $j=1,2,\ldots k$, we can go in one step from $x_1x_2\ldots x_i\ldots x_k$ to  $x_1x_2\ldots i\ldots x_k$, then the farthest vertices from $12\ldots k$ are of the form $\pi(1)\pi(2)\ldots \pi(k)$ for some permutation of $\{1,2,\ldots, k\}$ not fixing any element. Now consider the cyclic notation of $\pi$. Then, we claim that, for every cycle of length $\ell$, we need $\ell+1$ steps to put the corresponding digits in the `right' order.
For instance, in $F_{4\times 1}^1 (K_6)$ consider the vertex 4312=(1423). Thus, a path of length 5 to vertex 1234 is
$$
4312
\stackrel{(1)}{\sim} 5312 \stackrel{(4)}{\sim} 5314 
\stackrel{(2)}{\sim} 5214 
\stackrel{(3)}{\sim} 5234 
\stackrel{(1)}{\sim} 1234,
$$
where the superscript `$(i)$' means that we modify the $i$-th element.
Then, if $k$ is even, $k=2r$, a farthest vertex would be one with cyclic notation formed by $r$ cycles on length 2, which gives a path length of $3r$.
For instance, in $F_{4\times 1}^1 (K_6)$ a farthest vertex from 1234 will be 2143=(12)(34). Then a path of length 6 to vertex 1234 is
$$
2143
\stackrel{(1)}{\sim} 6143 \stackrel{(2)}{\sim} 6243 
\stackrel{(1)}{\sim} 1243 
\stackrel{(3)}{\sim} 1253 
\stackrel{(4)}{\sim} 1254
\stackrel{(3)}{\sim} 1234.
$$
Otherwise, if $\ell$ is odd, $\ell=2r+1$, a farthest vertex is one with a cyclic decomposition formed by $r-1$ cycles of length 2 plus one cycle of length 3. This gives a path length of $3(r-1)+4=3r+1$.
This completes the proof.
\end{proof}
 
\begin{figure}[t]
\begin{center}
\includegraphics*[width=10cm]{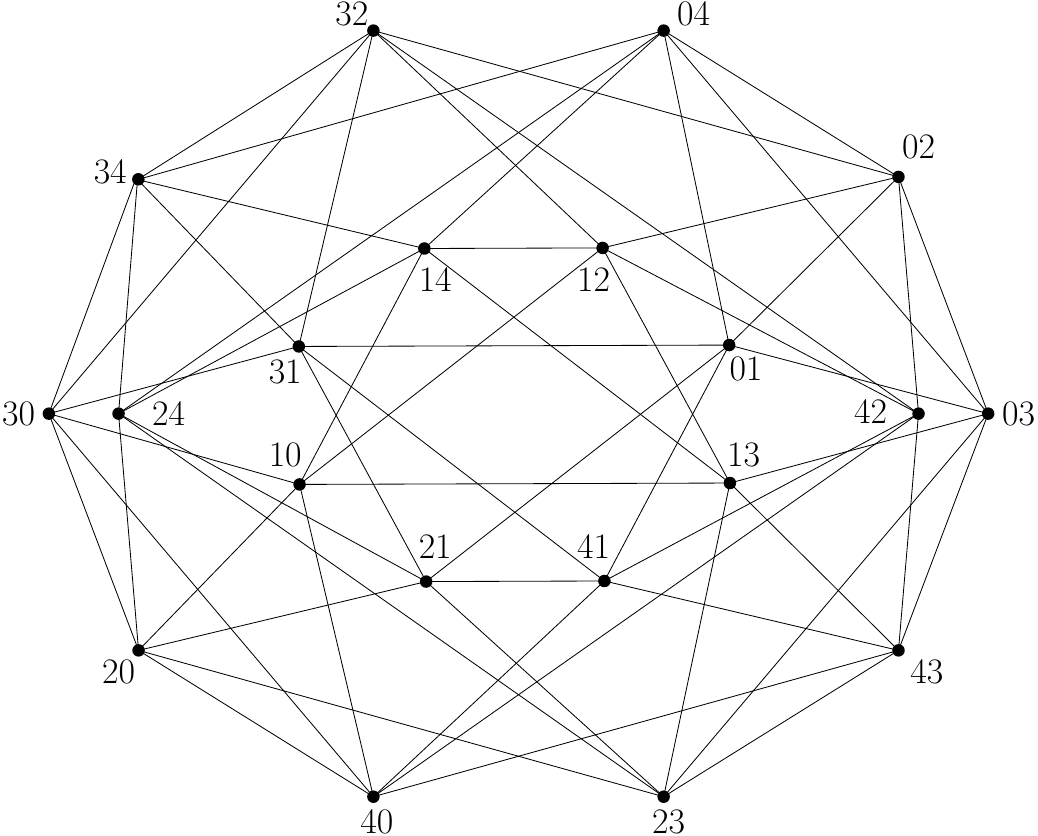}  
\caption{The supertoken graph $F_{2\times 1}^1(K_5)$ with diameter 3.} 
\label{fig:F_2x1^1(K5)}
\end{center}
\end{figure}

For instance, Figure \ref{fig:F_2x1^1(K5)} shows the supertoken graph $F_{2\times 1}^1(K_5)$ with diameter $3\frac{k}{2}=3$. Notice that the farthest vertex from $ij$ is $ji$ (in fact, this is an antipodal graph). Note also that, since $F_{k\times 1}^1(K_n)$ is vertex-transitive, its radius and diameter coincide.

\begin{figure}[t]
\begin{center}
\includegraphics*[width=6cm]{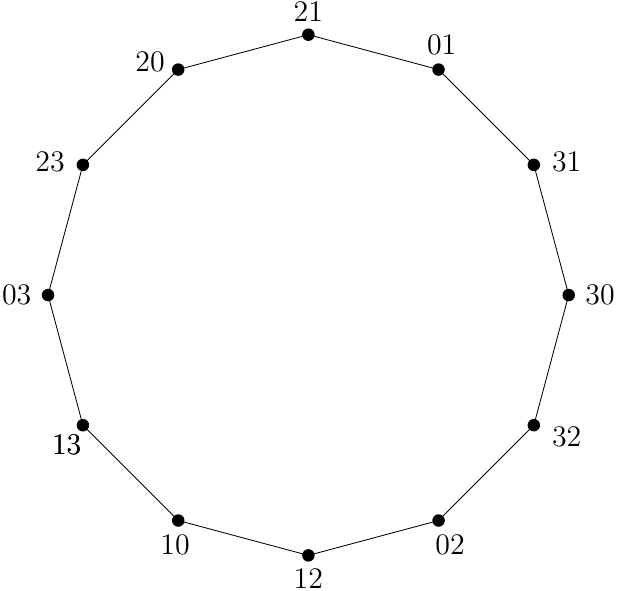}  
\caption{The supertoken graph $F_{2\times 1}^1(S_3)$ with diameter 6.} 
\label{fig:F_2x1^1(S3)}
\end{center}
\end{figure}

For the case of an $n$-pointed star $S_n$ ($n+1$ vertices) with central vertex $0$, we have a similar result. The proof is basically the same as the one of Proposition  \ref{propo:diamF2(Kn)}. The only difference is that the `right' movement of each token requires two steps instead of one. Thus, as shown in Figure \ref{fig:F_2x1^1(S3)} with the supertoken $F_{2\times 1}^1(S_3)$, to go from vertex $21$ (an involution) to $12$, we need 6 steps. As in the case of $F_{2\times 1}^1(K_n)$, the supertoken graph $F_{2\times 1}^1(S_n)$ is an antipodal graph with antipodal vertices $ij$ and $ji$, for $i,j\neq 0$.

\begin{proposition}
The diameter of the $k$-supertoken graph of $G=S_n$ is
$$
\diam (F_{k\times 1}^1(S_n))=
\left\{
\begin{array}{ll}
6r & \mbox{if $k=2r$,}\\
6r+1 & \mbox{if $k=2r+1$.\hskip 5cm  $\Box$}
\end{array}
\right.
$$
\end{proposition}

When the original graph is a cycle $C_n$ (see Figure \ref{fig:F_(2x1)^1(C_7)} for $F_{2\times 1}^1(C_7)$ and $F_{2\times 1}^1(C_6)$), we get the following result.

\begin{proposition}
The diameter of the 2-supertoken graph of the cycle $G=C_n$ is
$$
\diam(F_{2\times 1}^1(C_n))=n.
$$
\end{proposition}

\begin{figure}[t]
\begin{center}
\includegraphics*[width=10cm]{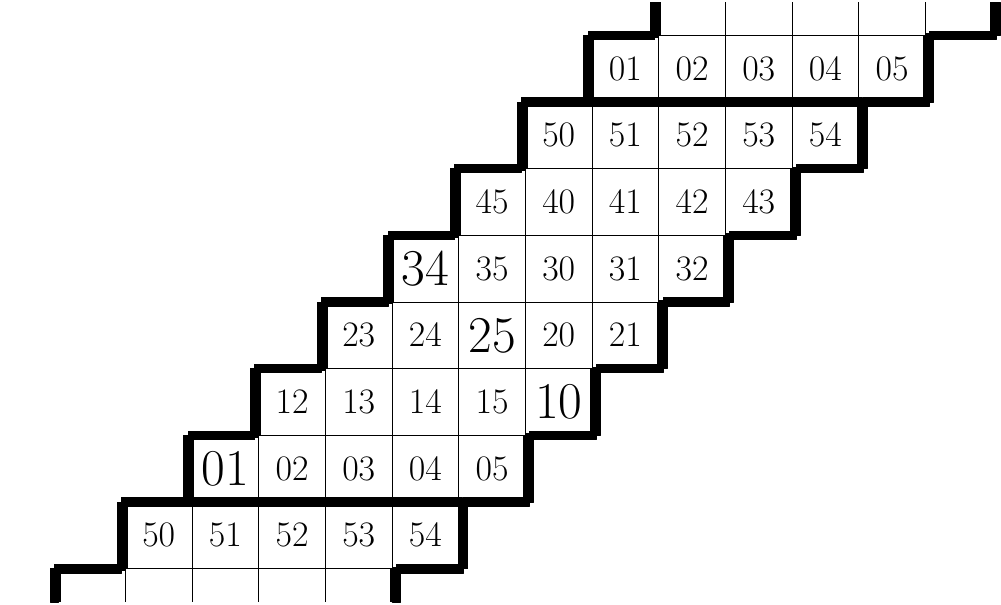}  
\caption{The farthest vertices from $01$ in $F_{2\times 1}^1(C_6)$ drawn as a plane tessellation are those in large-size labels.}
\label{fig:tessellation}
\end{center}
\end{figure}

\begin{proof}
 As before, a farthest vertex from $(i,j)$ is $(j,i)$. By symmetry, we can fix $i=0$. Then, we have the following paths of length $n$:
 $$
 \begin{array}{c}
 (0,1)\sim (0,2)\sim(0,3)\sim \cdots\sim (0,n-1)\sim(1,n-1)\sim(1,0);\\
(0,2)\sim(0,3)\sim\cdots\sim(0,n-1)\sim(1,n-1)\sim(2,n-1)\sim(2,0);\\
  \vdots \\
(0,n-1)\sim(1,n-1)\sim(2,n-1)\sim\cdots\sim(n-2,n-1)\sim(n-2,0)\sim(n-1,0).
\end{array}
 $$
 In fact, the farthest vertices from $(0,1)$, all of them at distance $n$, are $(1,0)$, $(2,n-1)$, $(3,n-2)$,\ldots,$(\lfloor \frac{n}{2}\rfloor, \lceil\frac{n}{2}\rceil+1)$. (See in Figure \ref{fig:tessellation} the farthest vertices from $01$ in $F_{2\times 1}^1(C_6)$ drawn as a plane tessellation). For more information about plane tessellations, see 
Yebra, Fiol, Morillo, and Alegre \cite{yfma85}.  
 This completes the proof.
\end{proof}

\section{Cayley graphs as supertoken graphs}
\label{sec:Cayley}
Recall that, given a group $\Gamma$ with generating set $\Delta$ such that $\Delta=\Delta^{-1}$, the Cayley graph $\Cay(\Gamma,\Delta)$ has vertices identified with the elements of $\Gamma$, and vertex $g$ is adjacent to the vertices $sg$ with $s\in \Delta$.
\begin{proposition}
\label{propo51}
    The supertoken graph $F_{(n-1)\times 1}^1 (K_n)$ is isomorphic to the bipartite Cayley graph on the symmetric group $S_n$ with generating set $\Delta=\{(1,n),(2,n),\ldots,(n-1,n)\}$.
\end{proposition}
\begin{proof}
  First, notice that $F_{(n-1)\times 1}^1 (K_n)$ has $n(n-1)\cdots 2=n!$ vertices as the number of elements of $S_n$.
  Moreover, each vertex $x_1x_2\ldots x_{n-1}$ can be identified with the permutation
  $x_1x_2\ldots x_{n-1}y$ where $y\neq x_1,\ldots,x_{n-1}$. Thus,
  $x_1x_2\ldots x_{n-1}y$ is adjacent to $yx_2\ldots x_{n-1}x_1$, $x_1y\ldots x_{n-1}x_2$, \ldots, $x_1x_2\ldots yx_{n-1}$. This corresponds to multiplying the permutation
  $x_1x_2\ldots x_{n-1}y$ by the involutions $(1,n)$, $(2,n)$, $\ldots$, $(n-1,n)$.
  This graph is bipartite with stable sets the even and odd permutations of $S_n$, since all the generators (involutions) are odd permutations.
\end{proof}
\begin{figure}[t]
\begin{center}
\includegraphics*[width=7cm]{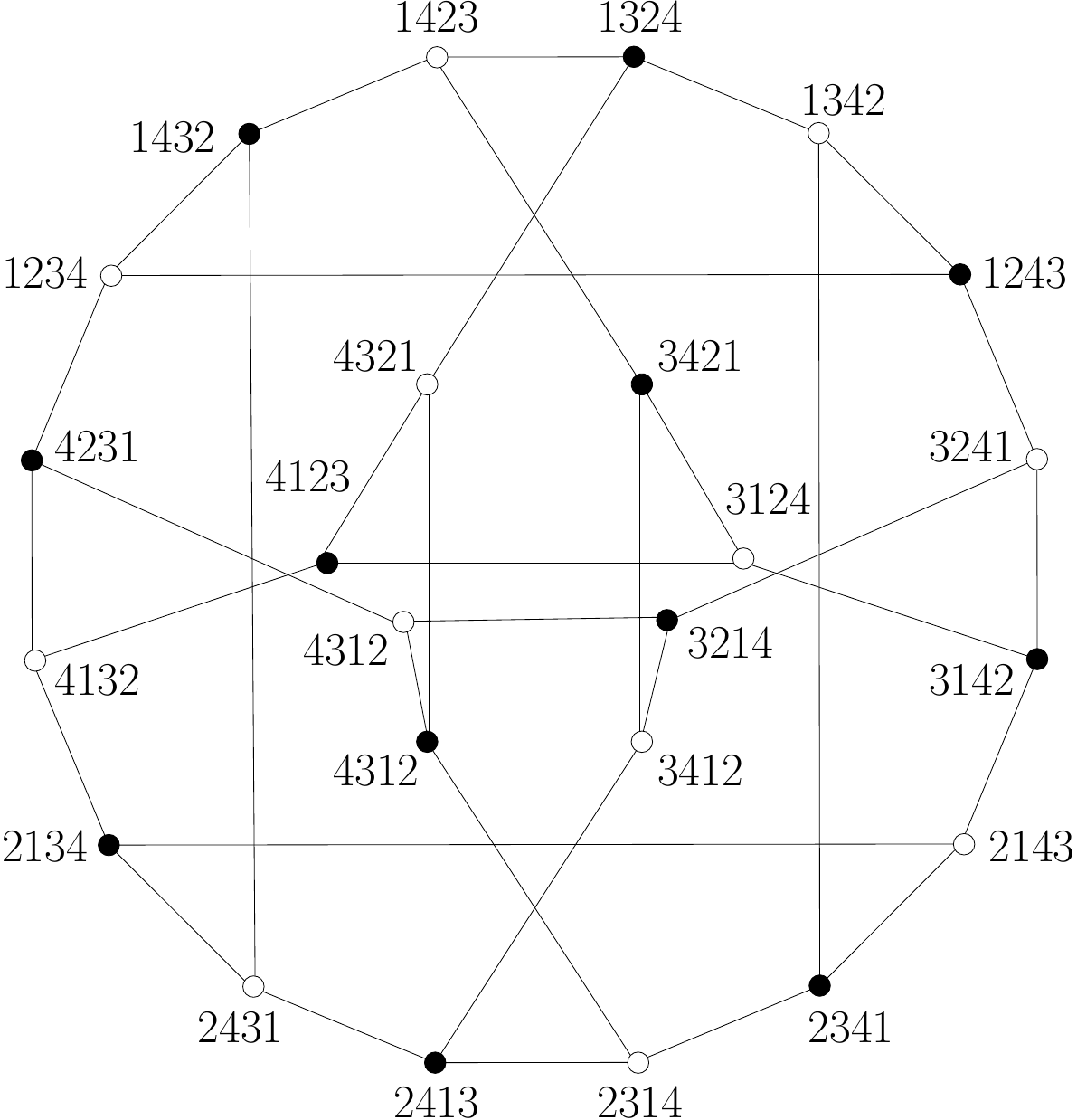}  
\caption{The supertoken (bipartite) graph $F_{3\times 1}^1(K_4)$ as a Cayley graph.} 
\label{fig:F_3x1^1(K4)}
\end{center}
\end{figure}
For instance, in Figure \ref{fig:F_3x1^1(K4)}, we show the supertoken graph $F_{3\times1}^1(K_4)$ as a Cayley graph on $S_4$.
From the above interpretation, we have the following consequence from  Proposition \ref{propo:diamF2(Kn)}.
\begin{corollary}
The diameter of the bipartite Cayley graph $\Cay(S_n,\Delta)$ with $\Delta=\{(1,n)$, $(2,n),\ldots,(n-1,n)\}$ is
\[
\diam\bigl(\Cay(S_n,\Delta)\bigr)
=\left\lfloor \frac{3n-3}{2} \right\rfloor
=
\begin{cases}
3\rho & \text{if } n=2\rho+1,\\
3\rho-2 & \text{if } n=2\rho.
\end{cases}
\]
\end{corollary}

\blue{
The Cayley graphs $\Cay(S_n,\Delta)$ are one of the examples used by 
 Araujo-Pardo, Del Río-Francos, López-Dudet,  Oliveros,  and Schulte \cite{adlos10},
to construct a vertex-transitive simple polytope, called the graphicahedron, whose 1-skeleton (edge graph) is $\Cay(S_n,\Delta)$.
}

\section{The spectra of supertoken graphs}
\label{sec:spectra}

Dalf\'o, Fiol, Mart\'inez-Barona, and Song studied the spectra of some supertoken graphs in \cite{dfms26}.
Thus, they proved the following two general results.

\begin{proposition}[\cite{dfms26}]
\label{prop:regularF1_1111}
    Let $G$ be a graph with order $n$. 
    Then, the vertex set of $F^{1}_{k\times 1}(G)$ admits a regular partition with quotient graph $F_k^1(G)$, the standard token graph.
\end{proposition} 

\begin{figure}[t]
\begin{center}
\includegraphics*[width=12cm]{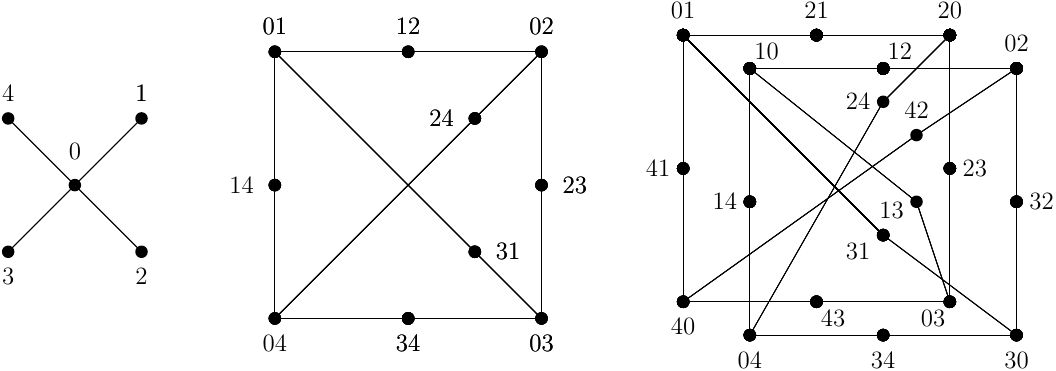}  
\caption{The star $S_4$, its 2-token graph $F_2(S_4)$, and its supertoken graph $F_{2\times 1}^1(S_4)$.} 
\label{fig:F_2(S4)}
\end{center}
\end{figure}

An example is shown in Figure \ref{fig:F_2(S4)}, 
where there is the star $S_4$, its 2-token graph $F_2(S_4)$, and its $2$-supertoken graph $F_{2\times 1}^1(S_4)$ with quotient graph $F_2(S_4)$.

As a consequence of this, the (adjacency or Laplacian) spectrum of $F_k^1(G)$ is contained in the corresponding spectrum of $F^{1}_{k\times 1}(G)$.

\begin{proposition}[\cite{dfms26}]
\label{propo:eigenvals F_{kx1}^1(G)}
    Let $\lambda_1,\ldots,\lambda_p$ be the (adjacency or Laplacian) eigenvalues of $F_k^1(G)$. Then, $\lambda_i$ is an (adjacency or Laplacian) eigenvalue of $F^{1}_{k\times 1}(G)$ for $i=1,\ldots,p=\binom{n}{k}$. 
    \end{proposition}

For instance, since $F_k^1(K_n)$ is isomorphic to the Johnson graph $J(n,k)$, we have the following corollary.
\begin{corollary}\label{coro63}
The spectrum of the supertoken graph $F_{k\times 1}^1(K_n)$ contains the eigenvalues of $J(n,k)$, namely,
$\lambda_j=(k-j)(n-k-j)-j$ with multiplicities $m_j=\binom{n}{j}-\binom{n}{j-1}$ for $j=0,1,\ldots,k$. In particular, $F_{k\times 1}^1(K_n)$ is a regular graph with degree $\lambda_0=k(n-k)$.
\end{corollary}

Now, we recall the following result by Verma \cite{v26}.
\begin{corollary}[\cite{v26}]
\label{coro:Verma}
    Let $G$ be a connected graph on $n$ vertices. Then, for any integer $k$ such that $2 \leq k \leq n$, the Laplacian spectral radius of the $k$-token graph $F_k(G)$ and the original graph $G$ coincide, that is, $\rho_L(F_k(G)) = \rho_L(G)$,
if and only if $G \simeq S_n$.
\end{corollary}

In the case of the Laplacian spectral radius of $F_{2\times 1}^1(G)$, we can say more, as follows.
\begin{proposition}
The Laplacian spectral radius of $F_{2\times 1}^1(G)$ coincides with the Laplacian spectral radius of $G$ \st{is}\greeen{if} and only if $G=S_n$, $\rho_L(F_{2\times 1}^1(S_n))=\rho_L(S_n)=n+1$.
\end{proposition}
\begin{proof}
This is a consequence of 
 Proposition \ref{propo:eigenvals F_{kx1}^1(G)}, and Corollary \ref{coro:Verma} by Verma \cite{v26}, who proved the same result for $F^{1}_{2}(S_n)$. From the Perron-Frobenius Theorem, the (adjacency or Laplacian) spectral radius of $F_k^1(G)$ and $F^{1}_{k\times 1}(G)$ coincide, as they both have positive eigenvectors.
 \end{proof}
 
Concerning the adjacency spectrum when the original graph is a star, Dalf\'o, Fiol, Mart\'inez-Barona, and Song \cite{dfms26} proved the following result.

\begin{proposition}[\cite{dfms26}]
\label{prop:sp-F_2x1^1(S_n)}
The adjacency spectrum of $F_{2\times 1}^1(S_n)$, where $S_n$ is the star with order $n+1$, is
$$ 
\spec F_{2\times 1}^1(S_n)=\left\{\pm \sqrt{2n-2}^{\,1},\pm \sqrt{n}^{\,n-1}, 
\pm \sqrt{n-2}^{\,n-1},
0^{\,n^2-3n+2}\right\}.
$$
\end{proposition}

In the following result, we provide the (adjacency) eigenvalues of the $k$-supertoken of a cycle graph.

\begin{figure}[t]
\begin{center}
\includegraphics*[width=8cm]{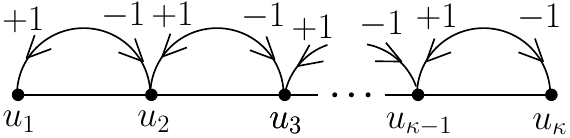}  
\caption{The base graph of the lift that gives $F_{k\times 1}^1(C_{m})$.} 
\label{fig:base-graph}
\end{center}
\end{figure}

\begin{proposition}
The complete spectrum of the supertoken graph $F_{k\times 1}^1(C_{m})$ is given by the eigenvalues
\begin{equation}
\label{eq:eigenvals}
\lambda(r,k)=4\cos\left(\frac{r\pi}{m}\right)\cos\left(\frac{k\pi}{m}\right), 
\end{equation}
for $r=0,\ldots,m-1$, and $k=1,\ldots,m-1$.
\end{proposition}
\begin{proof}
The supertoken graph $F_{k\times 1}^1(C_{m})$ can be seen as a lift on the group $\Z_{m}$ of the voltage graph on $\kappa=m-1$ vertices, as shown in Figure \ref{fig:base-graph}. The polynomial (tridiagonal) matrix $\BB(z)$, where $z=e^{r\frac{i2\pi}{m}}$, for $r=0,\ldots, m-1$, and its similar matrix $\BB^*(r)$
    are 
    $$
\BB(z) =
\left(
	\begin{array}{cccccc}
		0 & 1+z^{-1} & 0 & 0 & \ldots & 0 \\
		1+z & 0 & 1+z^{-1} & 0 & \ldots & 0 \\
		0 & 1+z & 0 & 1+z^{-1} & \ddots & 0 \\
		0 & 0 & 1+z & \ddots & \ddots  & 0 \\
		\vdots  & \vdots & \ddots &\ddots & 0 & 1+z^{-1}\\
		0 & 0 & \ldots &0 & 1+z & 0
\end{array}
\right)\cong
$$
	\begin{equation}
		\BB^*(r)=
		{\scriptstyle{\left(
		\begin{array}{cccccc}
			0 & 2\cos(\frac{r\pi}{m}) & 0 & 0 & \ldots & 0\\
			2\cos(\frac{r\pi}{m}) & 0 & 2\cos(\frac{r\pi}{m}) & 0 & \ldots & 0\\
			0 & 2\cos(\frac{r\pi}{m}) & 0 & 2\cos(\frac{r\pi}{m}) & \ddots& 0\\
			0 & 0 & 2\cos(\frac{r\pi}{m}) & \ddots &  \ddots & 0\\
			\vdots & \vdots & \ddots & \ddots
			& 0 & 2\cos(\frac{r\pi}{m}) \\[.1cm]
			0 & 0 & \ldots & 0 & 2\cos(\frac{r\pi}{m}) & 0
		\end{array}
		\right).}}
		\label{B-star}
	\end{equation}
For details about lift graphs and their spectra, see Dalf\'o, Fiol, Miller, Ryan, and \v{S}ir\'a\v{n} \cite{dfmrs17}, or Reyes, Dalf\'o, Fiol, and Messegu\'e \cite{rdfm23}. Now, by using the results of Losonczi \cite{l92} (rewritten later by Yueh \cite{y05}, as explained by Da Fonseca and Kowalenko \cite{fk20}), adjusted to our case, the eigenvalues of $\BB^*(r)$, and, hence, of the adjacency matrix of $F_{k\times 1}^1(C_{m})$, are as claimed.
\end{proof}
\begin{table}
\begin{tabular}{|c|cccccc|}
\hline
$r\backslash k$ & 1 & 2 & 3 & 4 & 5 & 6\\
\hline
0 & 3.60 & 2.49 & 0.89 & $-0.89$ & $-2.49$ & $-3.60$\\
1 & 3.24 & 2.24 & 0.80 & $-0.80$ & $-2.24$ & $-3.24$\\
2 & 2.24 & 1.55 & 0.55 & $-0.55$ & $-1.55$ & $-2.24$\\
3 & 0.80 & 0.55 & 0.19 & $-0.19$ & $-0.55$ & $-0.80$\\
4 & $-0.80$ & $-0.55$ & $-0.19$ & 0.19 & 0.55 & 0.80\\
5 & $-2.24$ & $-1.55$ & $-0.55$ & 0.55 & 1.55 & 2.24\\
6 & $-3.24$ & $-2.24$ & $-0.80$ & 0.80 & 2.24 & 3.24\\
\hline
\end{tabular}
\\
\caption{The eigenvalues of the supertoken graph $F^{1}_{2\times 1}(C_{7})$. }
\label{tab:eigenvalsC7}
\end{table}

\begin{table}
\begin{tabular}{|c|ccccc|}
\hline
$r\backslash k$ & 1 & 2 & 3 & 4 & 5\\
\hline
0 & 3.46 & 2 & 0 & $-2$ & $-3.46$\\
1 & 3 & 1.73 & 0 & $-1.73$ & $-3$\\
2 & 1.73 & 1 & 0 & $-1$ & $-1.73$ \\
3 & $0$ & 0 & 0 & 0 & 0 \\
4 & $-1.73$ & $-1$ & 0 & 1 & 1.73 \\
5 & $-3$ & $-1.73$ & 0 & 1.73 & 3 \\
\hline
\end{tabular}
\\
\caption{The eigenvalues of the supertoken graph $F^{1}_{2\times 1}(C_{6})$. }
\label{tab:eigenvalsC6}
\end{table}

For instance, the (adjacency) eigenvalues of the supertoken graphs $F^{1}_{2\times 1}(C_{7})$ and $F^{1}_{2\times 1}(C_{6})$ are shown in Tables \ref{tab:eigenvalsC7} and
\ref{tab:eigenvalsC6}, respectively.

\section{Independence, vertex Covers, and matchings in supertoken graphs}
\label{sec:ind-number}

\begin{figure}[t]
\begin{center}
\includegraphics*[width=14cm]{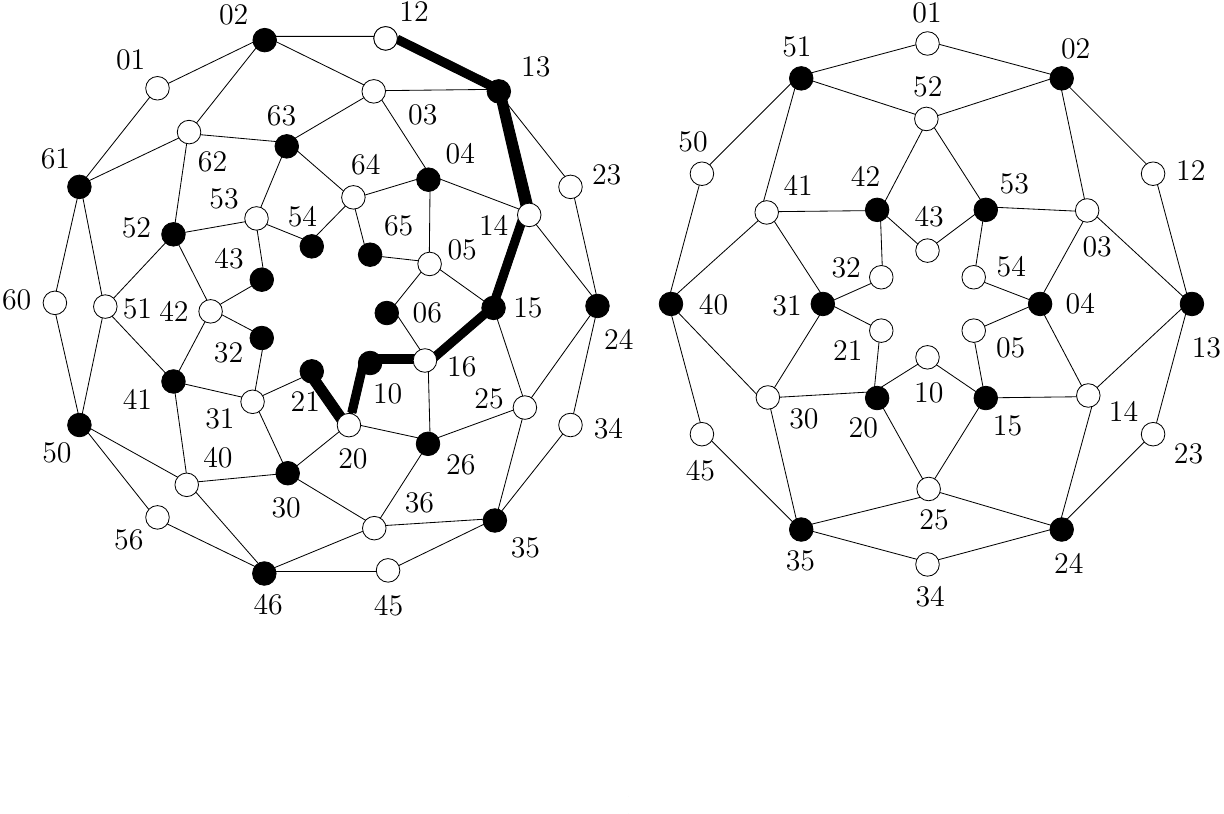}
\vskip-2cm
\caption{Some maximum independent sets (white vertices), of cardinalities $(2n+1)n=21$ and $2n^2=18$, together with some minimum vertex covers (black vertices), of cardinalities $(2n+1)n=21$ and $2n(n-1)=12$, in the supertoken graphs $F^{1}_{2\times 1}(C_{7})$ and $F^{1}_{2\times 1}(C_{6})$, respectively. {\purple{In a thick line, there is a path between the antipodal vertices $12$ and $21$.}} }
\label{fig:F_(2x1)^1(C_7)}
\end{center}
\end{figure}

In this section, we study several classical structural parameters of $k$-supertoken graphs, including the independence number, vertex covering number, and matching-related parameters. These parameters play an important role in graph theory, as they capture fundamental combinatorial properties of graphs and often reveal structural information about the original graph.

We begin by recalling the definitions of independent sets and vertex covers, together with their well-known relationship.

 A set of vertices in a graph \(G=(V,E)\) is called \emph{independent} if no two vertices in the set are adjacent. The \emph{independence number} of \(G\), denoted by \(\alpha(G)\), is the maximum cardinality of an independent set in \(G\). An independent set of cardinality \(\alpha(G)\) is called a \emph{maximum independent set}. Independent sets play a central role in graph theory and combinatorial optimization, as they represent collections of mutually non-adjacent vertices and are closely related to coloring problems and other structural properties of graphs.

A closely related concept is that of a vertex cover. A set \(C \subseteq V(G)\) is called a \emph{vertex cover} of \(G\) if every edge of \(G\) has at least one endpoint in \(C\). The \emph{vertex covering number} of \(G\), denoted by \(\beta(G)\), is the minimum cardinality of a vertex cover in \(G\). A vertex cover of cardinality \(\beta(G)\) is called a \emph{minimum vertex cover}.

These two notions are closely connected: the complement of any independent set is a vertex cover, and the complement of any vertex cover is an independent set. Consequently, if \(G\) has order \(n\), then
\[
\alpha(G)+\beta(G)=n.
\]
This identity highlights the natural duality between independent sets and vertex covers, and shows why both parameters are fundamental in the study of graph structure.

In \cite{aclr20}, de Alba, Carballosa, Leaños, and Rivera studied the
independence and matching numbers of some standard token graphs.
We now apply these notions to the study of $k$-supertoken graphs of bipartite graphs. In particular, we first determine the independence number of a bipartite $k$-supertoken graph $F^{1}_{k\times1}(G)$.

\begin{proposition}\label{prop51}
Let $G$ be a bipartite graph with stable sets $C_1$ and $C_2$, with cardinalities $|C_i|=c_i$ for $i=1,2$ such that $c_1\le c_2$. Then, given an integer $k(\ge 2)$,
the bipartite $k$-supertoken $F^{1}_{k\times 1}(G)$ has stable sets 
${\cal S}_1$ and ${\cal S}_2$.
Then, the independence number of $F^{1}_{k\times 1}(G)$ is
$$
\alpha(F^{1}_{k\times 1}(G))=\max\{|{\cal S}_1|,|{\cal S}_2|\},
$$
where
\begin{align}
|{\cal S}_1| &= \sum_{i=0}^{\min\{\lfloor c_1/2\rfloor,\lfloor k/2\rfloor\}}{k\choose 2i}\frac{c_1!}{(c_1-2i)!} \frac{c_2!}{(c_2-k+2i)!},\label{|S1|}\\
|{\cal S}_2| &= \sum_{i=0}^{\min\{\lfloor (c_1-1)/2\rfloor,\lfloor (k-1)/2\rfloor\}}{k\choose 2i+1}\frac{c_1!}{(c_1-2i-1)!} \frac{c_2!}{(c_2-k+2i+1)!}.
\label{|S2|}
\end{align}
\end{proposition}
\begin{proof}
Since the supertoken graph $F^{1}_{k\times 1}(G)$ is bipartite, its independence number equals the maximum cardinality of the stable sets, see Song, Dalf\'o, Fiol, Mora, and Zhang \cite{sdfmz26}. The cardinality of ${\cal S}_1$ in \eqref{|S1|} (respectively, ${\cal S}_2$ in \eqref{|S2|}) corresponds to the number of vertices having an even number $2i$ (respectively, odd number $2i+1$) of distinguished tokens in $C_1$. In the first case (the other being similar), for every $i$ in the range of the sum, we first choose the number ${k\choose 2i}$ of such tokens. Then, we compute the number of ways such tokens can be in $C_1$. That is, the variations of $c_1$ elements taken $2i$ at a time, $V_{c_1,2i}=\frac{c_1!}{(c_1-2i)!}$. Finally, we multiply the above numbers by the number of ways to choose the remaining tokens, $k-2i$, in $C_2$. That is, $V_{c_2,k-2i}=\frac{c_2!}{(c_2-(k-2i))!}$.  
\end{proof}

Recall that the $k$-supertoken graph $F^{1}_{k\times 1}(G)$ of a bipartite graph $G$ is also bipartite, but the converse is not necessarily true. See an example in Figure \ref{fig:F_(2x1)^1(C_7)} (left), in which $F_{2x1}^1(C_7)$ is bipartite, while $C_7$ is not.

The previous result provides an explicit expression for the independence number of $F^{1}_{k\times1}(G)$ when the underlying graph $G$ is bipartite. Since the complement of an independent set is a vertex cover, it is natural to ask how the vertex covering number behaves for these graphs.

\begin{proposition}
Let $G$ be a bipartite graph with stable sets $C_1$ and $C_2$, where $|C_i|=c_i$ for $i=1,2$ and $c_1\le c_2$. For every integer $k\ge 2$, if $F^{1}_{k\times 1}(G)$ is the bipartite $k$-supertoken graph of $G$ with stable sets ${\cal S}_1$ and ${\cal S}_2$, then
	\[
	\beta\!\left(F^{1}_{k\times 1}(G)\right)
	=
	\min\{|{\cal S}_1|,|{\cal S}_2|\}.
	\]
\end{proposition}

\begin{proof}
	Since $F^{1}_{k\times 1}(G)$ is bipartite with partite sets ${\cal S}_1$ and ${\cal S}_2$, we have
	\[
	|V(F^{1}_{k\times 1}(G))|=|{\cal S}_1|+|{\cal S}_2|.
	\]
	Moreover, by Proposition~\ref{prop51},
	\[
	\alpha\!\left(F^{1}_{k\times 1}(G)\right)=\max\{|{\cal S}_1|,|{\cal S}_2|\}.
	\]
	Recall that for any graph $H$,
	\[
	\alpha(H)+\beta(H)=|V(H)|.
	\]
	Applying this identity to $H=F^{1}_{k\times 1}(G)$ yields
	\[
	\beta\!\left(F^{1}_{k\times 1}(G)\right)
	=
	|V(F^{1}_{k\times 1}(G))|
	-
	\alpha\!\left(F^{1}_{k\times 1}(G)\right).
	\]
	Substituting the expressions above gives
	\[
	\beta\!\left(F^{1}_{k\times 1}(G)\right)
	=
	(|{\cal S}_1|+|{\cal S}_2|)-\max\{|{\cal S}_1|,|{\cal S}_2|\}
	=
	\min\{|{\cal S}_1|,|{\cal S}_2|\}.
	\]
\end{proof}

The previous proposition allows us to compute the independence and vertex covering numbers explicitly for several families of graphs. We illustrate this with some immediate consequences.

\begin{corollary}
		Let $G=S_n$ be the star graph of order $n+1$. Since $S_n$ is bipartite with partite sets of sizes $c_1=1$ and $c_2=n$, for every integer $k\ge 2$, the $k$-supertoken graph $F^{1}_{k\times 1}(G)$ satisfies
		\[
		\alpha\!\left(F^{1}_{k\times 1}(G)\right)
		=
		\max\left\{
		\frac{k\,c_2!}{(c_2-k+1)!},
		\frac{c_2!}{(c_2-k)!}
		\right\},
		\]
		and
		\[
		\beta\!\left(F^{1}_{k\times 1}(G)\right)
		=
		\min\left\{
		\frac{k\,c_2!}{(c_2-k+1)!},
		\frac{c_2!}{(c_2-k)!}
		\right\}.
		\]
\end{corollary}

\begin{corollary}
	\label{coro:2-super-bip}
		Let $G$ be a bipartite graph with stable sets $C_1$ and $C_2$, where $|C_i|=c_i$ for $i=1,2$ and $c_1\le c_2$. Then, for the $2$-supertoken graph $F^{1}_{2\times 1}(G)$, we have
		\[
		\alpha\!\left(F^{1}_{2\times 1}(G)\right)
		=
		\max\{2c_1c_2,\; c_1(c_1-1)+c_2(c_2-1)\},
		\]
		and
		\[
		\beta\!\left(F^{1}_{2\times 1}(G)\right)
		=
		\min\{2c_1c_2,\; c_1(c_1-1)+c_2(c_2-1)\}.
		\]
\end{corollary}

We next specialize to the case where the original graph is a cycle and determine both the independence number and the vertex covering number of the corresponding $2$-supertoken graph.

\begin{proposition}
Let $G=C_m$ be the cycle graph of order $m$. Then the independence number of its $2$-supertoken graph $F^{1}_{2\times 1}(C_m)$ is
\[
\alpha(F^{1}_{2\times 1}(C_m))=
\begin{cases}
	(2n+1)n, & \text{if } m=2n+1,\\
	2n^2, & \text{if } m=2n,
\end{cases}
\]
and, consequently, the vertex covering number is
\[
\beta(F^{1}_{2\times 1}(C_m))=
\begin{cases}
	(2n+1)n, & \text{if } m=2n+1,\\
	2n(n-1), & \text{if } m=2n.
\end{cases}
\]
\end{proposition}

\begin{proof}
	First observe that the order of the graph $F^{1}_{2\times 1}(C_m)$ depends on the parity of $m$. If $m=2n$, then 
	\[
	|V(F^{1}_{2\times 1}(C_m))|=2n(2n-1),
	\]
	while if $m=2n+1$, then
	\[
	|V(F^{1}_{2\times 1}(C_m))|=2n(2n+1).
	\]
	We first determine the independence number of the corresponding $2$-supertoken graph, distinguishing two cases depending on the parity of $m$.
	
	\medskip
	
	\noindent
	\textbf{Even case.} Suppose that $m=2n$. Then $C_{2n}$ is bipartite with partite sets of equal size $n$. Hence, the result follows directly from Corollary~\ref{coro:2-super-bip} with $c_1=c_2=n$. In particular, a maximum independent set is formed by the vertices
	\[
	(i,i+1),\,(i,i+3),\,\ldots,\,(i,i+2n-1), \qquad i\in\mathbb{Z}_{2n}.
	\]
	An illustration of this construction is shown in Figure~\ref{fig:F_(2x1)^1(C_7)} (right) for the case $C_6$ ($n=3$).
	
	\medskip
	
	\noindent
	\textbf{Odd case.} Suppose that $m=2n+1$. In this case, an independent set of $F^{1}_{2\times 1}(C_{2n+1})$ is again given by the vertices
	\[
	(i,i+1),\,(i,i+3),\,\ldots,\,(i,i+2n-1), \qquad i\in\mathbb{Z}_{2n+1}.
	\]
	See Figure~\ref{fig:F_(2x1)^1(C_7)} (left) for the example corresponding to $C_7$ ($n=3$).
	
	To show that this independent set is maximum, note that the graph $F^{1}_{2\times 1}(C_{2n+1})$ has $N=(2n+1)2n$ vertices. Moreover, by \eqref{eq:eigenvals}, its spectrum contains $N/2$ positive eigenvalues and $N/2$ negative eigenvalues. Therefore, by the Cvetkovi\'c bound~\cite{c71},
	\[
	\alpha\le \frac{N}{2}=(2n+1)n.
	\]
	Since the independent set constructed above has cardinality $(2n+1)n$, it follows that it is maximum.
	
	\medskip
	
	Finally, recall that for every graph $H$,
	\[
	\alpha(H)+\beta(H)=|V(H)|.
	\]
	Applying this identity to $H=F^{1}_{2\times 1}(C_m)$ gives
	\[
	\beta\!\left(F^{1}_{2\times 1}(C_m)\right)
	=
	|V(F^{1}_{2\times 1}(C_m))|
	-
	\alpha\!\left(F^{1}_{2\times 1}(C_m)\right),
	\]
	from which the statement follows by substituting the expressions obtained above.
\end{proof}

We now turn to another fundamental family of graphs, namely paths, and determine the corresponding parameters for their $2$-supertoken graphs.

\begin{corollary}
	Let $G=P_m$ be the path graph on $m\ge 2$ vertices. Then
	\[
	\alpha\!\left(F^{1}_{2\times 1}(P_m)\right)
	=
	\begin{cases}
		\dfrac{m^2}{2}, & \text{if } m \text{ is even},\\[2mm]
		\dfrac{m^2-1}{2}, & \text{if } m \text{ is odd},
	\end{cases}
	\]
	and
	\[
	\beta\!\left(F^{1}_{2\times 1}(P_m)\right)
	=
	\begin{cases}
		\dfrac{m(m-2)}{2}, & \text{if } m \text{ is even},\\[2mm]
		\dfrac{(m-1)^2}{2}, & \text{if } m \text{ is odd}.
	\end{cases}
	\]
\end{corollary}

\begin{proof}
This is a consequence of Corollary \ref{coro:2-super-bip}.
\end{proof}

Beyond independence and vertex covers, it is also natural to study matching-related parameters of supertoken graphs. We therefore recall some standard definitions and classical results concerning matchings and edge covers.

A \emph{matching} in a graph $G=(V,E)$ is a set of pairwise non-adjacent edges; that is, no two edges in the set share a common endpoint. A matching is said to be \emph{maximum} if it has the largest cardinality among all matchings of $G$. The cardinality of a maximum matching is called the \emph{matching number} of $G$ and is denoted by $\nu(G)$.

Equivalently, a matching may be regarded as an independent set of edges. In this context, the \emph{edge independence number} of $G$, denoted by $\alpha'(G)$, is the maximum cardinality of a set of pairwise non-adjacent edges. Thus,
\[
\alpha'(G)=\nu(G).
\]

A closely related notion is that of an \emph{edge cover}. An edge cover of a graph $G$ is a set of edges such that every vertex of $G$ is incident with at least one edge of the set. The \emph{edge covering number} of $G$, denoted by $\beta'(G)$, is the minimum cardinality of an edge cover of $G$. An edge cover of cardinality $\beta'(G)$ is called a \emph{minimum edge cover}.

These two parameters are related by a classical identity: if $G$ is a graph of order $n$ with no isolated vertices, then
\[
\alpha'(G)+\beta'(G)=n.
\]
This equality expresses the natural duality between maximum matchings and minimum edge covers.

In the case of bipartite graphs, matchings are also closely related to vertex covers. Indeed, the K\H{o}nig--Egerv\'ary theorem 
(see Deming \cite{d79} and Kung \cite{k84})
states that, for every bipartite graph $G$,
\[
\nu(G)=\beta(G),
\]
where $\beta(G)$ denotes the minimum cardinality of a vertex cover of $G$. Since $\alpha'(G)=\nu(G)$, it follows that
\[
\alpha'(G)=\nu(G)=\beta(G).
\]

Moreover, for every graph $G$,
\[
\alpha(G)+\beta(G)=|V(G)|.
\]
Therefore, if $G$ is bipartite and has no isolated vertices, then
\[
\beta'(G)=|V(G)|-\nu(G)=|V(G)|-\beta(G)=\alpha(G).
\]
Hence, in a bipartite graph without isolated vertices, the edge covering number coincides with the independence number. In particular, the size of a maximum independent set is completely determined by the size of a maximum matching.

These classical relations apply naturally to $k$-supertoken graphs whenever the original graph is bipartite. In particular, when $G$ is bipartite, they yield useful connections among matchings, vertex covers, and edge covers of $F^{1}_{k\times1}(G)$, as shown below.



\begin{proposition}
\label{prop:KE-supertoken}
Let $G$ be a bipartite graph. Then
\[
\alpha'\!\left(F^{1}_{k\times1}(G)\right)=\beta\!\left(F^{1}_{k\times1}(G)\right).
\]
Moreover, if $F^{1}_{k\times1}(G)$ has no isolated vertices, then
\[
\beta'\!\left(F^{1}_{k\times1}(G)\right)
=
\alpha\!\left(F^{1}_{k\times1}(G)\right).
\]
\end{proposition}

\begin{proof}
    The statement follows directly from the preceding observations and the K\H{o}nig--Egerv\'ary theorem.
\end{proof}

\section{{Clique number of supertoken graphs}}
\label{sec:clique}

The \emph{clique number} of a graph $G$, denoted by $\omega(G)$, is the maximum
size of a complete subgraph of $G$. Determining or bounding the clique number
is a classical problem in graph theory and has been studied for several graph
transformations, including token graphs.

For the standard $k$-token graph $F_k(G)$, Fabila-Monroy, Flores-Pe\~{n}aloza,
Huemer, Hurtado, Urrutia, and Wood~\cite{ffhhuw12} determined the exact value
of the clique number.

\begin{theorem}[\cite{ffhhuw12}]
    $\omega(F_k(G)) = \min\{\omega(G), \max\{n - k + 1, k + 1\}\}.$
\end{theorem}
For $k$-supertoken graphs, we first obtain a relationship between the clique numbers of the supertoken and token graphs. Then, we provide the exact value of the clique number of a supertoken graph. 
\begin{proposition}
\label{thm:clique}
Let $G$ be a graph on $n$ vertices and $k \ge 1$. 
%
Then,
$$
\omega\big(F^{1}_{k\times 1}(G)\big) \leq \omega\big(F_k(G)\big),
$$
where $F_k(G)$ is the standard $k$-token graph of $G$ (with vertices as sets instead of vectors). 
\end{proposition}
\begin{proof}
Define the projection
$$
\pi : V\big(F^{1}_{k\times 1}(G)\big) \to V\big(F_k(G)\big),
\qquad
\pi(v_1,\dots,v_k) := \{v_1,\dots,v_k\}.
$$
Let $X$ be a clique in $F^{1}_{k\times 1}(G)$. For any two vertices 
$\vu=(u_1,\dots,u_k)$ and $\vv=(v_1,\dots,v_k)$ in $X$, adjacency means that $\vv$ is obtained from $\vu$ by replacing exactly one coordinate $u_i$ by a vertex $v_i\in V(G)$ such that $u_iv_i \in E(G)$, and all other coordinates remain unchanged.\\
Hence,
$$
\pi(\vv) = \big(\pi(\vu)\setminus\{u_i\}\big)\cup\{v_i\},
$$
with $u_iv_i \in E(G)$. Therefore, $\pi(\vu)$ and $\pi(\vv)$ are adjacent in $F_k(G)$.
Thus, $\pi(X)$ is a clique in $F_k(G)$, and consequently
$$
|X|=|\pi(X)| \le \omega(F_k(G)).
$$
\end{proof}



\begin{proposition}\label{prop83}
The clique number of the supertoken graph is
\begin{equation}
\label{eq:clique}
\omega\big(F^{1}_{k\times 1}(G)\big)=\min\{\omega(G), n - k + 1\}. 
\end{equation}
\end{proposition}
\begin{proof}
First, we prove that in $F^{1}_{k\times 1}(G)$ there is only one type of cliques. Namely, if
$Y$ is a clique in $F_{k\times 1}^1(G)$, we have that $Y = \{S \cup \{v\} : v \in K\}$, where $S \subseteq V(G)$ with $|S|=k-1$ and $K$ is a clique in $G$ disjoint from $S$. This follows from the fact that, in this type of supertokens, to form a clique, only one token may vary while the remaining $k-1$ tokens stay fixed on the same vertices.
Then, \eqref{eq:clique} holds from the two following cases:
\begin{itemize}
\item 
If $k\le n-\omega(G)+1$, then $\omega\big(F^{1}_{k\times 1}(G)\big)=\omega(G)$. In this case, $S$ can be chosen such that $S\cap K=\emptyset$, with $K$ being a clique in $G$ of size $\omega(G)$.
\item 
If $k= n-c+1$ where $c<\omega(G)$, then $\omega\big(F^{1}_{k\times 1}(G)\big)=c=n-k+1$. Now, $S$ can be chosen such that
$S\cap K=\emptyset$, with $K$ being a clique in $G$ of size $c$.
\end{itemize}
\end{proof}

\section{{Chromatic number of supertoken graphs}}
\label{sec:chrom}

The \emph{chromatic number} of a graph $G$, denoted by $\chi(G)$, is the minimum
number of colors needed to color the vertices of $G$ so that adjacent vertices
receive different colors. Studying how the chromatic number behaves under graph transformations is a classical topic in graph theory. 

For token graphs, Fabila-Monroy, Flores-Peñaloza, Huemer, Hurtado, Urrutia, and
Wood~\cite{ffhhuw12} obtained the following bounds.

\begin{theorem}[\cite{ffhhuw12}]\label{t81}
 Let $G$ be a graph on $n$ vertices. Then,
 $$
 \left(\frac{1}{2}+\frac{2}{n}\right)\chi(G)-1\le \chi(F_k(G))\le \chi(G)
 $$
 for every $k\ge 1$.
 \end{theorem}
Concerning the supertoken graph, we have the following result.
\begin{proposition}\label{prop62}
Let $G$ be a graph of order $n$ and let $k \ge 1$. Then,
$$
\omega\!\left(F^{1}_{k\times 1}(G)\right) \leq \chi(F_{k\times 1}^1(G))\le
\chi(F_k(G)) \le \chi(G).
$$
\end{proposition}

\begin{proof}
The inequality $\omega\!\left(F^{1}_{k\times 1}(G)\right) \leq \chi(F_{k\times 1}^1(G))$
follows from the general fact that, for any graph, the chromatic number is at least its clique number. To show that $\chi(F_{k\times 1}^1(G))\le
\chi(F_k(G))$, we observe that
every vertex-coloring of $F_k(G)$ can be easily extended to a coloring of $F_{k\times 1}^1(G)$. Indeed, in any supertoken graph $F_{k\times 1}^1(G)$, we have $k!$ times the number of vertices of a token graph $F_k(G)$. That is, given a vertex $v_1v_2\ldots v_k$ in $F_k(G)$, we have the $k!$ permutations of $v_1v_2\ldots v_k$ in $F_{k\times 1}^1(G)$, corresponding to the $k!$ vertices for any vertex in $F_k(G)$. Moreover, these $k!$ vertices ($v_1v_2\ldots v_k$ included) can be colored with the same color as $v_1v_2\ldots v_k$, as they are never adjacent. So, $\chi(F_{k\times 1}^1(G))\le \chi(F_k(G))$. Finally, the inequality $\chi(F_k(G))\le \chi(G))$ follows from Theorem~\ref{t81}.
\end{proof}

\begin{corollary}
	Let \(G\) be a graph of order \(n\), and let \(F^{1}_{k\times 1}(G)\) be its \(k\)-supertoken graph. Then
	\[
	\min \bigl\{\omega(G),\, n-k+1\bigr\}
	\;\le\;
	\chi\!\left(F^{1}_{k\times 1}(G)\right)
	\;\le\;
	\chi(G).
	\]
\end{corollary}

\begin{proof}
    This is an immediate consequence of Propositions~\ref{prop83} and~\ref{prop62}.
\end{proof}

\section{Planarity of Supertoken Graphs}
\label{sec:planar}
A graph is said to be \emph{planar} if it admits an embedding in the plane in which no two edges intersect except possibly at a common endpoint. Such an embedding is called a \emph{plane embedding} of the graph. Recall that every planar graph with \(n\ge 3\) vertices and \(m\) edges satisfies $m \le 3n-6$. Moreover, if the graph is bipartite, then the stronger inequality $m \le 2n-4$
holds.

\begin{lemma}\label{leq5}
Let \(G\) be a planar graph that is regular of degree \(\ell\). Then \(\ell \le 5\).
\end{lemma}

\begin{proof}
Let \(n\) and \(m\) denote the number of vertices and edges of \(G\), respectively. Since \(G\) is \(\ell\)-regular, we have
\[
2m = \ell n.
\]
On the other hand, since \(G\) is planar and \(n \ge 3\), it follows that
\[
m \le 3n - 6.
\]
Combining these inequalities yields
\[
\ell n = 2m \le 6n - 12,
\]
and hence
\[
\ell \le 6 - \frac{12}{n} < 6.
\]
Therefore, \(\ell \le 5\).
\end{proof}

A fundamental characterization of planar graphs is provided by Kuratowski's theorem.

\begin{theorem}[Kuratowski's Theorem]
A graph is planar if and only if it contains no subgraph which is a subdivision of \(K_5\) or \(K_{3,3}\).
\end{theorem}

Here, a \emph{subdivision} of a graph is obtained by replacing its edges with internally vertex-disjoint paths.

These classical results will be used throughout this section to determine whether the supertoken graphs under consideration are planar or nonplanar.

We start by analyzing the supertoken graphs \(F^{1}_{2\times 1}(K_n)\) for \(n \ge 2\).

\begin{proposition}\label{prop104}
  The graph \(F^{1}_{2\times 1}(K_n)\) is planar if and only if \(n \in \{2,3,4\}\).
\end{proposition}

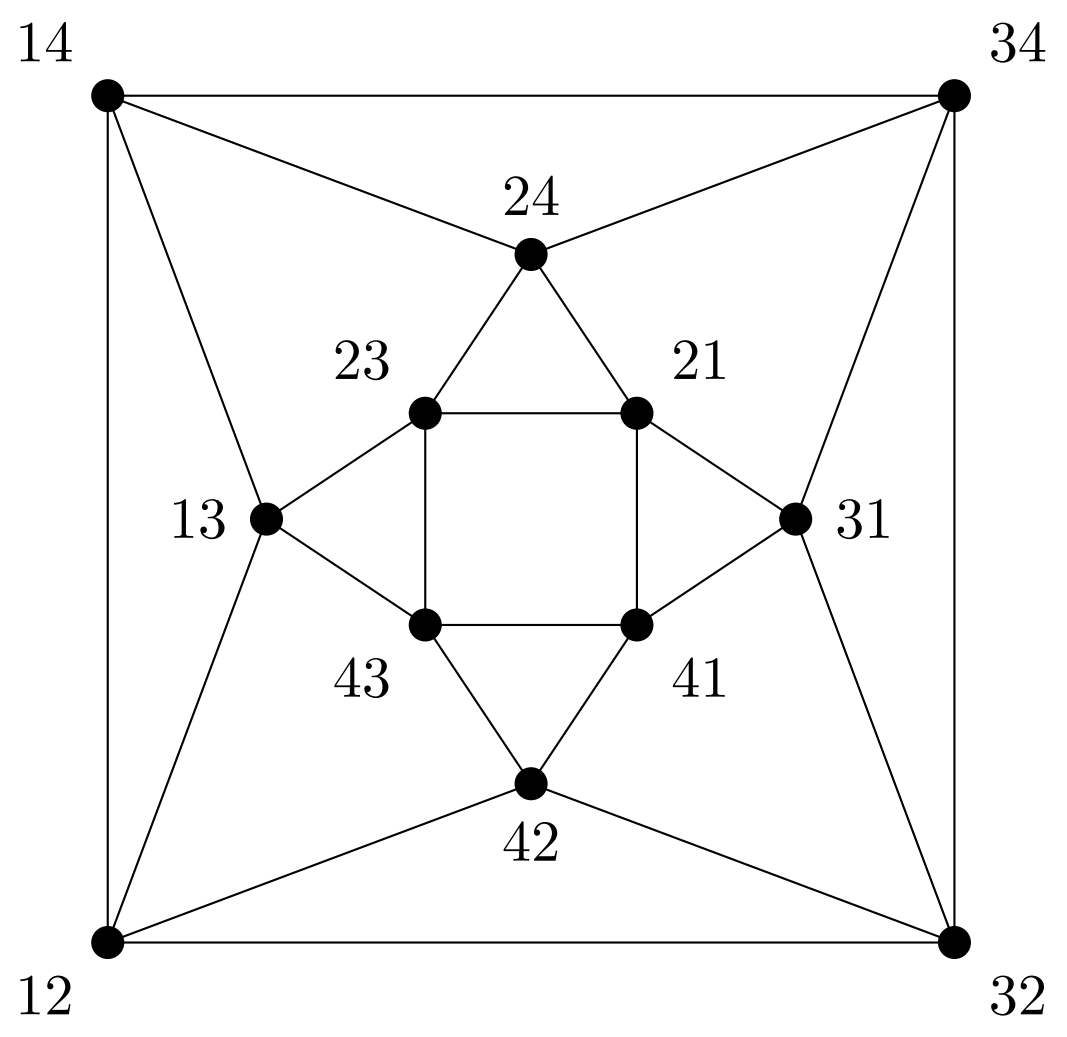
\begin{figure}[ht]
\centering
\begin{tikzpicture}[
    scale=0.7,
    every node/.style={font=\small},
    v/.style={circle, fill=black, inner sep=2.2pt}
]
\coordinate (14) at (-4,4);
\coordinate (34) at ( 4,4);
\coordinate (12) at (-4,-4);
\coordinate (32) at ( 4,-4);
\coordinate (24) at (0, 2.5);
\coordinate (13) at (-2.5,0);
\coordinate (31) at ( 2.5,0);
\coordinate (42) at (0,-2.5);
\coordinate (23) at (-1,1);
\coordinate (21) at ( 1,1);
\coordinate (43) at (-1,-1);
\coordinate (41) at ( 1,-1);
\draw (14)--(34)--(32)--(12)--cycle;
\draw (14)--(24)--(34);
\draw (14)--(13)--(12);
\draw (34)--(31)--(32);
\draw (12)--(42)--(32);
\draw (24)--(23) (24)--(21);
\draw (13)--(23) (13)--(43);
\draw (31)--(21) (31)--(41);
\draw (42)--(43) (42)--(41);
\draw (23)--(21)--(41)--(43)--cycle;
\node[v,label=above left:$14$] at (14) {};
\node[v,label=above right:$34$] at (34) {};
\node[v,label=below left:$12$] at (12) {};
\node[v,label=below right:$32$] at (32) {};

\node[v,label=above:$24$] at (24) {};
\node[v,label=left:$13$] at (13) {};
\node[v,label=right:$31$] at (31) {};
\node[v,label=below:$42$] at (42) {};

\node[v,label=above left:$23$] at (23) {};
\node[v,label=above right:$21$] at (21) {};
\node[v,label=below left:$43$] at (43) {};
\node[v,label=below right:$41$] at (41) {};

\end{tikzpicture}
\caption{A plane embedding of 
\(F^{1}_{2\times 1}(K_4)\). }
\label{fig:F21K4}
\end{figure}

\begin{proof}
First, suppose that the graph \(F^{1}_{2\times 1}(K_n)\) is planar for some \(n\ge 2\). 
By Corollary~\ref{coro63}, the graph \(F^{1}_{2\times 1}(K_n)\) is regular of degree 
\(\ell = 2(n-2)\). By Lemma~\ref{leq5}, every planar regular graph has degree at most \(5\), 
and therefore \(\ell \le 5\). It follows that
\[
2(n-2)\le 5,
\]
which implies \(n \le \tfrac{9}{2}\). Since \(n\) is an integer, we conclude that 
\(n\in \{2,3,4\}\).

Conversely, if \(n\in \{2,3,4\}\), then \(F^{1}_{2\times 1}(K_n)\) is planar. Indeed,
\(F^{1}_{2\times 1}(K_2)\cong 2K_1\), the complement of \(K_2\), and 
\(F^{1}_{2\times 1}(K_3)\cong C_6\), the cycle on six vertices; both graphs are planar. 
For the case \(n=4\), see Figure~\ref{fig:F21K4}, where a plane embedding of 
\(F^{1}_{2\times 1}(K_4)\) is exhibited. 
\end{proof}

Let $G$ be a graph on $n$ vertices. By Remark~\ref{prop:any-graph}, when \(k=n\), the graph \(F^{1}_{k\times 1}(G)\) consists of \(n!\) isolated vertices and is therefore planar. We now turn to the case \(1\le k\le n-1\), with \(n\ge 2\), and study the planarity of \(F^{1}_{k\times 1}(G)\) in this range.

\begin{lemma}\label{lem3}
The graph \(F^{1}_{3\times 1}(K_4)\) is non-planar.
\end{lemma}

\begin{proof}
The vertices of \(F^{1}_{3\times 1}(K_4)\) are the ordered triples \((a,b,c)\) of distinct elements of \(\{1,2,3,4\}\), where two vertices are adjacent if and only if one coordinate is replaced by the unique element not appearing in the triple.

Consider the sets
\[
A=\{(3,1,4),\ (3,4,1),\ (4,2,1)\}
\quad\text{and}\quad
B=\{(2,3,1),\ (3,1,2),\ (3,2,1)\}.
\]
We show that these vertices serve as the branch vertices of a subdivision of \(K_{3,3}\).

gFor every pair \(u\in A\) and \(v\in B\), there exists a path in \(F^{1}_{3\times 1}(K_4)\) joining \(u\) to \(v\), and these paths can be chosen to be internally vertex-disjoint. For instance, the following paths realize all such connections:
\[
(3,1,4)-(2,1,4)-(2,3,4)-(2,3,1), \quad
(3,1,4)-(3,1,2), \quad
(3,1,4)-(3,2,4)-(3,2,1),
\]
\[
(3,4,1)-(2,4,1)-(2,3,1), \quad
(3,4,1)-(3,4,2)-(3,1,2), \quad
(3,4,1)-(3,2,1),
\]
\[
(4,2,1)-(4,3,1)-(2,3,1), \quad
(4,2,1)-(4,2,3)-(4,1,3)-(4,1,2)-(3,1,2), \quad
(4,2,1)-(3,2,1).
\]
\begin{figure}[t]
\begin{center}
\includegraphics*[width=8cm]{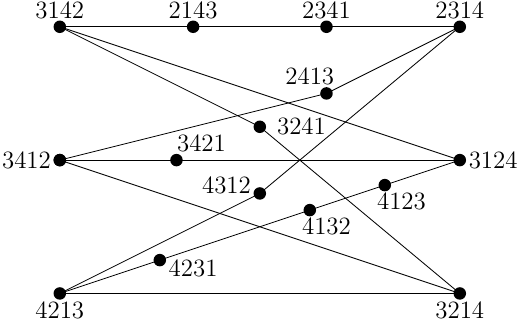} 
\caption{A subdivision of $K_{3,3}$ in 
\(F^{1}_{3\times 1}(K_4)\), where we added the final digit in every vertex so that it can represent a permutation (see the proof of Proposition \ref{propo51} in Section \ref{sec:Cayley}). }
\label{fig:K33}
\end{center}
\end{figure}
(See Figure \ref{fig:K33}.)
Hence \(F^{1}_{3\times 1}(K_4)\) contains a subdivision of \(K_{3,3}\). By Kuratowski's theorem, it follows that \(F^{1}_{3\times 1}(K_4)\) is non-planar.
\end{proof}

\begin{theorem}
Let \(n\ge 2\) and \(1\le k\le n-1\). The graph \(F^{1}_{k\times 1}(K_n)\) is planar if and only if
\[
(n,k)\in\{(2,1),(3,1),(3,2),(4,1),(4,2)\}.
\]
\end{theorem}

\begin{proof}
Assume that \(n\ge 2\) and \(1\le k\le n-1\), and suppose that
\(F^{1}_{k\times 1}(K_n)\) is planar. By Corollary~\ref{coro63},
the graph \(F^{1}_{k\times 1}(K_n)\) is regular of degree
\[
\ell = k(n-k).
\]
Hence, by Lemma~\ref{leq5}, we must have \(\ell \le 5\), and therefore
\[
k(n-k)\le 5.
\]
It follows that $\{(n,k): n\ge 2,\ 1\le k\le n-1,\ k(n-k)\le 5\}$ is
\[\{(2,1),(3,1),(3,2),(4,1),(4,2),(4,3),(5,1),(5,4),(6,1),(6,5)\}.
\]

We now examine these possibilities.

For \(k=1\), we have \(F^{1}_{1\times 1}(K_n)\cong K_n\). Since \(K_5\) is non-planar and occurs as a subgraph of \(K_n\) for all \(n\ge 5\), it follows that \(F^{1}_{1\times 1}(K_n)\) is planar if and only if \(n\le 4\). Thus the admissible pairs are \((2,1),(3,1),(4,1)\).

For \(k=2\), Proposition~\ref{prop104} implies that the only candidates are \(n=3\) and \(n=4\), yielding the pairs \((3,2)\) and \((4,2)\).

By Lemma~\ref{lem3}, the pair \((4,3)\) is excluded.

It remains to consider the pairs \((5,4)\) and \((6,5)\). By Proposition~\ref{prop51}, the graph \(F^{1}_{k\times1}(K_n)\) is bipartite in both cases. Let \(N(n,k)\) and \(M(n,k)\) denote its number of vertices and edges, respectively. Since \(F^{1}_{k\times1}(K_n)\) is a \(k(n-k)\)-regular graph, we have
\(
M(n,k)=\frac{N(n,k)}{2}\,k(n-k).
\)
As \(k(n-k)\ge 4\), it follows that
\[
M(n,k)>2N(n,k)-4.
\]

Hence, \(F^{1}_{k\times1}(K_n)\) exceeds the planar bipartite bound and is therefore non-planar. In particular, this rules out the cases \((5,4)\) and \((6,5)\).

Consequently, the only planar cases are
\[
(2,1),(3,1),(3,2),(4,1),(4,2).
\]
A direct inspection shows that each of these graphs is planar, completing the proof.
\end{proof}


 


\subsection*{Acknowledgments} 
This research work has been funded by AGAUR, the Catalan Government, under project 2021SGR00434, and by MICINN, the Spanish Government, under project PID2020-115442RB-I00.
The first and fourth authors were also supported by a grant from the Research Vice-Rector of Universitat de Lleida. In addition, this work was partially supported by the Slovenian Research and Innovation Agency (ARIS) under research program P1-0285 and research projects J1-50000 and Z1-70008.

\subsection*{Conflict of interest}
The authors declare no conflict of interest. 

\subsection*{Data availability statement} 
All data generated or analyzed during this study are included in this article.



\begin{thebibliography}{99}
	
\bibitem{adlos10}
G. Araujo-Pardo, M. Del Río-Francos, M. L\'opez-Dudet, D. Oliveros, E. Schulte,
The graphicahedron, {\em European J. Combin.} {\bf 31} (2010) 1868--1879.

\bibitem{agrr07}
K. Audenaert, C. Godsil, G. Royle, and T. Rudolph,
Symmetric squares of graphs,
\textit{J. Combin. Theory B} \textbf{97} (2007) 74--90.

\bibitem{bdfs24}
E. T. Baskoro, C. Dalf\'o, M. A. Fiol, and R. Simanjuntak,
On some metric properties of supertoken graphs, 
{\textit{Util. Math.} \textbf{123} (2024) 19--33.}

\bibitem{c71}
D. M. Cvetkovi\'c, 
Graphs and their spectra, 
\textit{Publ. Elektrotehn. Fak. Ser. Mat. Fiz.} {\bf 354-356} (1971) 1--50.

\bibitem{fk20}
C. M. Da Fonseca and V. Kowalenko,
Eigenpairs of a family of tridiagonal matrices: three decades later, 
\textit{Acta Math. Hungar.} \textbf{160} (2020) 376--389.

\bibitem{dfms26}
C. Dalf\'o, M. A. Fiol, B. Mart\'inez-Barona, and X. Song,
On the spectra of generalized token graphs of stars,
Submitted, 2026.
 
\bibitem{dfmrs17}
C. Dalf\'o, M. A. Fiol, M. Miller, J. Ryan, and J. \v{S}ir\'a\v{n},
An algebraic approach to lifts of digraphs,
\textit{Discrete Appl. Math.} {\bf 269} (2019) 68--76.

\bibitem{aclr20}
H. de Alba, W. Carballosa, J. Leaños, and L. M. Rivera,
Independence and matching numbers of some token graphs, 
\textit{Australas. J. Combin.} {\bf 76(3)} (2020) 387--403.

\bibitem{d79}
R. W. Deming, 
Independence Numbers of Graphs--An Extension of the Koenig-Egervary Theorem,
\textit{Discrete Math.} \textbf{27} (1979) 23--33.

\bibitem{ffhhuw12}
R. Fabila-Monroy, D. Flores-Pe\~{n}aloza, C. Huemer, F. Hurtado, J. Urrutia, and D. R. Wood,
Token graphs,
\textit{Graphs Combin.} \textbf{28(3)} (2012) 365--380.



\bibitem{HaSm17}
R. H. Hammack and G. D. Smith,
Cycle bases of reduced powers of graphs,
\textit{Ars Math. Contemp.} \textbf{12} (2017) 183--203.

\bibitem{k84}
J. P. S. Kung, 
Jacobi's Identity and the König-Egerváry Theorem,
\textit{Discrete Math.} \textbf{49} (1984) 75--77.

\bibitem{l92}
{L. Losonczi, 
Eigenvalues and eigenvectors of some tridiagonal matrices, 
\textit{Acta Math. Hung.} {\bf 60} (1992) 309--332.}



\bibitem{rdfm23} 
M. A. Reyes, C. Dalf\'o, M. A. Fiol, and A. Messegu\'e,
On the spectra of token graphs of cycles and other graphs,
\textit{Linear Algebra Appl.} \textbf{679} (2023) 38--66.

\bibitem{sdfmz26}
{X. Song, C. Dalfó, M. A. Fiol, M. Mora, and S. Zhang, 
On generalized token graphs,
\textit{Filomat} {\textbf{40:2} (2026) 721--738}.}

\bibitem{v26} 
P. Verma, Laplacian spectral radius and integrality of token graphs, {\em Electron. J. Linear Algebra} \textbf{42} (2026) 27--38.


\bibitem{yfma85}
J. L. A. Yebra, M. A. Fiol, P. Morillo, and I. Alegre,
The diameter of undirected graphs associated to plane tessellations,
\textit{Ars Combin.} {\bf 20B} (1985) 159--171.

\bibitem{y05}
{W.-C. Yueh, 
Eigenvalues of several tridiagonal matrices, 
\textit{Appl. Math. E-Notes} \textbf{5} (2005) 66--74.}

\end{thebibliography}
\end{document}